\documentclass[10pt,a4paper,twoside]{article}
\usepackage{amssymb}
\usepackage{latexsym}
\usepackage{fullpage}
\usepackage{a4,latexsym,rotating,amsmath,epsfig}
\usepackage{here}
\usepackage{amsmath}
\usepackage{amscd}
\usepackage{oldgerm}
\usepackage{fancyheadings}
\usepackage{here}

\def\R{\mathbb{R}}

\def\E{\mathfrak{E}}

\def\SO{\mathbf{SO}\,}

\def\SU{\mathbf{SU}\,}

\def\Sp{\mathbf{Sp}\,}
\def\G2{\mathbf{G_2}\,}
\def\Spin{\mathbf{Spin}\,}
\def\1{\mathbf{1}}
\def\:{\lrcorner}
\def\#{\sharp}
\def\k{\kappa}

\def\a{\alpha}
\def\b{\beta}

\def\d{\delta}

\def\s{\sigma}
\def\t{\tau}

\def\<#1,#2>{\langle#1,\,#2\rangle}

\def\L{\mathcal{L}}

\def\E{\mathcal{E}}

\def\v{\noindent}

\def\h{\mathfrak{h}}
\def\curv{\mathcal{R}}
\def\ric{Ric}
\def\s{\mathbf{\t}}
\def\W{\mathcal{W}}
\def\dsum{\oplus}

\def\tens{\otimes}
\def\Hol{Hol}
\def\sumONB{\sum_{\scriptscriptstyle\mathrm{ONB}}}

\def\qed{\ensuremath{\quad\Box\quad}}

\def\inv#1{\raise.1em\hbox to 0pt{$^{-1}$\hss}_{#1}\;}

\newtheorem{Theorem}{Theorem}
\newtheorem{Lemma}[Theorem]{Lemma}
\newtheorem{Corollary}[Theorem]{Corollary}

\newtheorem{Example}[Theorem]{Example}
\newtheorem{Remark}[Theorem]{Remark}
\newtheorem{Proposition}[Theorem]{Proposition}

\title{The geometry of conformally Einstein metrics with degenerate Weyl tensor}\author{Jesse Alt}

\begin{document}

\maketitle

\begin{center}

Humboldt-Universit\"at zu Berlin, Germany, alt@math.hu-berlin.de \\ Supported by the DFG International Research Training Group GRK 870

\end{center}

\begin{abstract}
The problem of characterizing conformally Einstein manifolds by tensorial conditions has been tackled recently in papers by M. Listing, and in work by A. R. Gover and P. Nurowski. Their results apply to metrics satisfying a ''non-degeneracy" condition on the Weyl tensor $\W$. We investigate the geometry of the foliations arising on conformally Einstein spaces (with Riemannian signature) where this condition fails, which then allows us to characterize a general class of locally conformally Einstein Riemannian manifolds with degenerate Weyl tensor.
\end{abstract}

\section{Introduction}
The importance of conformal geometry for physics and mathematical physics is seen, for example, in the AdS/CFT correspondence. On the other hand, Einstein spaces are obviously of interest for these fields, as well as for differential geometry. It's thus a natural problem to try to characterize (Riemannian) manifolds which are conformally Einstein. That is: given a Riemannian manifold $(M^n,g)$, does there exist a smooth function $\phi$ on $M$ such that the conformally related metric $\bar{g} = e^{2\phi}g$ is Einstein? Criteria for answering this question locally were investigated for arbitrary dimension $n$ by Listing in \cite{lis01}, and the main result was later generalized by Gover and Nurowski in \cite{GovNur} and by Listing in \cite{lis04}. These results require the following "non-degeneracy" condition on the Weyl tensor $\W$: If we define $Ker(\W)_p = \{X \in T_p(M) : X \: \W = 0\}$, then $Ker(\W)_p$ must have rank $0$ on a dense subset of $M$.\\

We say a Riemannian manifold $(M,g)$ has {\it rank $k$ Weyl degeneracy} if $Ker(\W)_p$ has rank $k$ on a dense subset of $M$. It follows from this definition that there is an open, dense subset $M_k \subset M$, on which $Ker(\W)$ is a distribution. If $(M,g)$ is conformally Einstein, then this distribution has nice geometric properties. It is a foliation (integrable), totally umbilical - and with respect to an Einstein metric in the conformal class, it is even totally geodesic. Moreover, if $k \geq 4$, then the leaves of the foliation are conformally flat. These properties follow from the identification (see section 4), for an Einstein metric, of $Ker(\W)$ with the so-called $\k$-nullity distribution. This allows us to draw conclusions for the conformal geometry. The central result, established in sections 6 and 7, is:

\begin{Theorem} \label{main} Let $(M,g)$ be a locally conformally Einstein Riemannian manifold with rank $k$ Weyl degeneracy. Then, locally, $(M_k,g)$ falls into one of two classes - transversally integrable or transversally non-integrable. In the transversally integrable case, if the Einstein scaling is non-negative, then $(M_k,g)$ must be locally conformally Ricci-flat and have decomposable conformal holonomy. In the transversally non-integrable case, there is a unique Einstein metric in the conformal class, and the conformal Einstein re-scaling can be determined locally by an algorithm.
\end{Theorem}

In section 2, we explain where the non-degeneracy restriction on $Ker(\W)$ comes from, and review the main concepts in Listing's method. There are several reasons why it's interesting to study conformally Einstein spaces where this restriction does not hold. The first comes from the conformal tractor calculus, facts from which are briefly reviewed in section 3. For an Einstein metric, the equation $X \: \W = 0$ is the integrability condition for the existence of parallel tractors. $Ker(\W)$ can be seen as parametrizing the space of possible normal conformal Killing 1-forms on a conformally Einstein space (cf. \cite{lei04}). The next, related motivation is the example of Einstein-Sasaki spaces $(M^{2n+1},g,\xi)$, where the Reeb vector $\xi$ of the Sasaki contact structure is seen to satisfy $\xi \: \W = 0$. Einstein-Sasaki spaces, which are of course fascinating from the purely differential geometric viewpoint, have also been intensely studied in mathematical physics (cf. \cite{mar06}). We review in section 4 some basic facts from contact and Sasakian geometry, and give in section 5 a tensorial characterization of conformally Einstein-Sasaki spaces with rank 1 Weyl degeneracy. Further, although we limit ourselves here to Riemannian signature, it should be noted that spaces with degenerate Weyl tensor arise in Lorentzian geometry from Brinkmann waves (spaces with a parallel light-like vector field $V$). For pp waves, Fefferman waves and for all Brinkmann waves which have twistor-spinors, the parallel lightlike vector field satisfies $V \: \W = 0$. The simplest examples of Brinkmann waves, plane waves, are known to be conformally Ricci-flat (cf. \cite{leis05}).\\

\section{Listing's method for characterizing conformally Einstein metrics}

To get conditions for a Riemannian manifold to be (locally) conformally Einstein, Listing constructs a tensor $\mathbb{T}$  which is a candidate to be the (local) gradient of a positive function giving a conformal transformation to an Einstein metric. The obstacle to applying this method to metrics with degenerate Weyl tensor is that the candidate $\mathbb{T}$ is no longer unique. We review here the main ideas of Listing's proof.\\

On a Riemannian manifold $(M,g)$ we have the usual invariants associated to the metric: the Levi-Civita connection $\nabla$, Riemannian curvature tensor $\curv$, as well as Ricci curvature $\ric$ and scalar curvature $\s$. We write $\nabla f$ for the gradient of a smooth function $f$. When we want to emphasize or clarify the dependence on a particular choice of metric, we add sub- or superscripts, e.g. $\nabla^g, \curv^g, \s_g$, etc. In addition, we have:\\

\v the {\it traceless Ricci tensor}
\begin{align*}
\ric^o := \ric - (n-1) \s g
\end{align*}

\v the {\it Schouten tensor}
\begin{align*}
\h := \frac {1} {n-2} (\ric - \frac {\s} {2(n-1)}g);
\end{align*}

\v the {\it Kulkarni-Nomizu product} of two symmetric (2,0) tensors $A$ and $B$
\begin{align*}
A \star B (X,Y,U,V) := A(X,U)&B(Y,V) + A(Y,V)B(X,U)\\
&- A(X,V)B(Y,U) - A(Y,U)B(X,V);
\end{align*}

\v the {\it Weyl tensor}
\begin{align*}
\W := \curv - g \star \h;
\end{align*}

\v the {\it Cotton-York tensor}
\begin{align*}
C(X,Y,Z) := (\nabla_X \h)(Y,Z) - (\nabla_Y \h)(X,Z);
\end{align*}

\v the {\it divergence} of the Weyl tensor is given by
\begin{align*}
\d \W (X,Y,Z) := \sumONB (\nabla_{E_i} \W) (X,Y,Z, E_i);
\end{align*}

\v and the {\it rotation} of a vector field $V$ (we write $V^*$ for the dual 1-form to $V$) is the 2-form
\begin{align*}
rot(V)(X,Y) := dV^*(X,Y) = g(\nabla_X V,Y) - g(X,\nabla_Y V).
\end{align*}

A space with $\d \W = 0$ is called a {\it C-space}, or space with {\it harmonic Weyl tensor}. Noting the identity, $\d \W = (n-3)C$, we see that, in particular, Einstein spaces are C-spaces (cf. \cite{Besse}, p. 440, for a discussion of the sense in which C-spaces generalize Einstein spaces). The key to Listing's method is finding a uniquely determined tensor which must vanish if $(M,g)$ is conformally equivalent to a C-space, which is an integrability condition for the conformal Einstein condition. This relies on the behaviour of $\d \W$ under conformal translation. Namely, for a conformally related metric $\bar{g} = e^{2\phi}g$, we have:
\begin{align}\label{conf cspace}
\bar{\d} \bar{\W} = \d \W + (3-n) (\nabla^g \phi) \: \W.
\end{align}

Thus an integrability condition for $(M,g)$ to be locally conformally Einstein is: there must be a local gradient field $X$ (i.e., $rot(X)=0$) such that 

\begin{align*}
\d \W + (n-3) X \: \W = 0.
\end{align*}

\v Note that if the homogeneous part of this equation ($(n-3)X\:\W=0$) has only the trivial solution, then a solution to these integrability conditions, if it exists, must be unique. This is where the non-degeneracy condition on the Weyl tensor comes from. When it is fulfilled, Listing constructs a vector field $\mathbb{T}$, which must be the unique solution, if one exists, to (\ref{conf cspace}).\\ 

Define the tensor $F_V$, determined by a vector field $V$ as follows:
\begin{align*}
F_V(X,Y) := g(\nabla_X V,Y) + g(X,V)g(Y,V) - \frac {1} {n} [div(V) + g(V,V)]g(X,Y).
\end{align*}
\v From the transformation behavior of $\ric$ under a conformal change of metric, it follows that a metric $e^{2\phi}g$ is Einstein if and only if $\ric_g^o + (2-n) F_{\nabla^g \phi} = 0$. On the other hand, a vector field $V$ is a local gradient if and only if $F_V$ is symmetric. Using these facts, Listing concludes that a space with non-degenerate Weyl tensor is locally conformally Einstein if and only if $\ric_g^o + (2-n) F_{\mathbb{T}} = 0$.\\

\section{The tractor calculus for conformal geometry and conformally Einstein metrics}

We note that an improved result by Listing in \cite{lis04} was established by Gover and Nurowski in \cite{GovNur} using the so-called standard tractor calculus for conformal structures. We summarize here briefly the essential properties of the tractor calculus, following \cite{lei04}, in order to formulate the important results of F. Leitner and S. Armstrong on the conformal holonomy of conformally Einstein metrics. For more details on the background to the constructions needed for the tractor calculus see, e.g. \cite{CSS} or \cite{fehl}.\\

The tractor calculus for conformal geometry relies on the construction of a vector bundle along with connection, $(\mathcal{T}(M), \nabla^{NC})$, both of which are canonically determined by the class of conformally equivalent metrics $c = [g]$ on $M$. The fibers of $\mathcal{T}$ are isomorphic to $\R^{(p+1,q+1)}$, where $(p,q)$ is the signature of $c$, and comes equipped with an inner product. In the standard manner for vector bundles, we can construct the dual of the tractor bundle, and exterior products, to get bundles of $p$-forms on tractors, denoted $\Omega^p_{\mathcal{T}}(M)$. The canonical connection, $\nabla^{NC}$, is also extended to these bundles in the standard way. The construction gives important conformal invariants derived from $\nabla^{NC}$ (for instance, the curvature of this connection, $\curv^{\nabla^{NC}}$, and the conformal holonomy, $\Hol^{\nabla^{NC}} = Hol(M,c)\subset\SO (p+1,q+1)$).\\

Fixing a metric $g \in c$, we get a direct-sum decomposition of the forms on tractors in the usual differential forms on $M$:
\begin{align*}
\Omega^{p+1}_{\mathcal{T}}(M) = \Omega^p(M) \dsum \Omega^{p+1}(M) \dsum \Omega^{p-1}(M) \dsum \Omega^p(M)
\end{align*}
\v Writing $\nabla^{NC}$ in terms of the decomposition given by a metric $g$, we get the following, nice form for the action of covariant differentiation on tractor-forms:
$$\nabla^{NC}_X =
\left(\begin{array}{cccc} 
  \nabla^{LC}_X & -X\: & X^*\wedge & 0 \\
  -\h(X)^*\wedge & \nabla^{LC}_X & 0 & X^* \wedge \\
  \h (X) \: & 0 & \nabla^{LC}_X & X\: \\
  0 & \h (X) \: & \h(X)^* \wedge & \nabla^{LC}_X
\end{array}\right)$$

We are interested in forms which are parallel with respect to $\nabla^{NC}$, and Leitner shows that such forms are uniquely determined by the leading component under the decomposition of $\Omega^{p+1}_{\mathcal{T}}(M)$ given by a choice of metric in the conformal class. This $p$-form is called a {\it normal conformal Killing $p$-form}, and Leitner gives a set of four equations characterizing these forms. The existence of normal conformal Killing forms, which are conformally covariant, gives important information about the conformal geometry. Note that the curvature of the normal conformal connection, $\curv^{\nabla^{NC}}$, also has a nice form which gives further integrability conditions for the existence of normal conformal Killing forms:
$$\curv^{\nabla^{NC}}(X,Y) =
\left(\begin{array}{cccc} 
  \W(X,Y) & 0 & 0 & 0 \\
  -C(X,Y)^* \wedge & \W(X,Y) & 0 & 0 \\
  C(X,Y) \: & 0 & \W(X,Y) & 0 \\
  0 & C(X,Y) \: & C(X,Y)^* \wedge & \W(X,Y)
\end{array}\right)$$

A smooth function $\a$ on $M$ which corresponds under $g$ to the first component of a $\nabla^{NC}$-parallel 1-form on tractors (a normal conformal Killing function), gives, away from its zero set, a conformal transformation to an Einstein metric. Since a non-trivial function satisfying the normal conformal Killing equations is non-zero almost everywhere, we have a following correspondence between vectors in $\R^{(p+1,q+1)}$ fixed by $\Hol^{\nabla^{NC}}$ and Einstein structures in the conformal class $c$ which are defined up to singularities of measure zero on $M$. Furthermore, the causality of the $\Hol^{\nabla^{NC}}$-fixed vector determines the scalar curvature of the corresponding Einstein metric in the conformal class: for spacelike vectors, the scalar curvature is negative; for lightlike vectors, Ricci-flat; and for time-like vectors, scalar curvature is positive. The structure determined by a parallel tractor is thus called an {\it almost Einstein metric}.\\

Using independent methods, Leitner and Armstrong each established the following fact, which makes it possible to introduce a notion of decomposability for conformal holonomy:

\begin{Theorem}
(Armstrong \cite{arm05}, Leitner \cite{lei04}) Let $(M,c)$ be a conformal manifold with Riemannian signature. There exists a $\Hol (M,c)$-invariant subspace of $\R^{(n+1,1)}$ of dimension $p$, $2 \leq p \leq n$, if and only if there exists a metric $g \in c$ such that:\\

\v 1. $(M,g)$ is locally isometric to $(H^p \times L^q, h + l)$, with $n=p+q$ and $(H,h)$, $(L,l)$ almost Einstein;\\

\v 2. The scalar curvatures are related by: $\s_h = \frac {-p(p-1)} {q(q-1)} \s_l$;\\

\v 3. The conformal holonomies are related by: $\Hol(M,[g]) \cong \Hol(H,[h]) \times \Hol(L,[l])$.
\end{Theorem}

A conformal manifold for which these conditions hold is called {\it (conformally) decomposable} (or: it has decomposable conformal holonomy). An {\it indecomposable} conformal Riemannian manifold, therefore, has conformal holonomy which fixes no subspace of (co-)dimension greater than one. Considering conformally Einstein manifolds with Riemannian signature that are indecomposable, Armstrong went on to give the following classification results for conformal holonomy:

\begin{Theorem}
(Armstrong \cite{arm05}) Let $(M,g)$ be an indecomposable conformally Einstein manifold of Riemannian signature, with non-zero scalar curvature in the Einstein scaling. Then the conformal holonomy is one of the following:\\
\indent -$\SO(n,1)$ \\
\indent -$\SO(n+1)$ \\
\indent -$\SU(\frac {n+1} {2})$ \\
\indent -$\Sp(\frac {n+1} {4}) $ \\
\indent -$\G2 $ \\
\indent -$\Spin (7) $ \\

In case $(M,g)$ is conformally Ricci-flat, then the conformal holonomy is isomorphic to a semi-direct product of the Riemannian holonomy with $\R^n$, $Hol(M,[g]) \cong Hol(M,g) \rtimes \R^n$, and the possible indecomposable conformal holonomy groups are:\\
\indent -$\SO(n) \rtimes \R^n$ \\
\indent -$\SU(\frac {n} {2}) \rtimes \R^n$ \\
\indent -$\Sp(\frac {n} {4}) \rtimes \R^n$ \\
\indent -$\G2 \rtimes \R^7$ \\
\indent -$\Spin (7) \rtimes \R^8$ \\

\v Moreover, all these holonomy groups do actually occur.\\
\end{Theorem}

We note briefly the significance of Einstein-Sasaki manifolds in this picture. In \cite{sem}, Semmelmann showed how a Sasaki structure defines a special type of (conformal) Killing form. Namely, he defines a Sasaki structure on a Riemannian manifold $(M,g)$ to be given by a unit length Killing vector field $\xi$ satisfying for any vector field $X$ the equation
\begin{align*}
\nabla_X ( d\xi^* ) = -2 X^* \wedge \xi^*,
\end{align*}
\v and shows that this is equivalent to the classic definition of Sasaki structure from contact geometry (see below). The contact form $\eta$ which is dual to $\xi$ is a Killing 1-form by definition, while the additional condition says that $\eta$ is a special Killing form for the constant $-2$. Semmelmann proves a one-one correspondence between special Killing $p$-forms on the Riemannian manifold $(M,g)$ and parallel $(p+1)$-forms on the metric cone $(\hat{M},\hat{g})$. In \cite{arm05}, on the other hand, the conformal holonomy classification for Einstein spaces with non-zero Einstein scaling, is shown by constructing an isomorphism between the standard tractor bundle associated to such manifolds, and the metric cone. From this it follows that conformally Einstein-Sasaki manifolds give the first example of special conformal holonomy $\SU(\frac {n+1} {2})$, since this is a well-known characteristic of the holonomy of the metric cone of an Einstein-Sasaki space.\\

The same can be seen at the level of normal conformal Killing forms. As noted by Leitner, the normal conformal Killing equations reduce, for a co-closed $p$-form $\b$ on an Einstein space, to
\begin{align*}
\nabla_X \b &= \frac {1} {p+1} X \: d\b \\
\nabla_X d\b &= -\frac {(p+1)\s} {n(n-1)} X^* \wedge \b.
\end{align*}
In our case, the contact $1$-form $\eta$ of the Sasaki structure is by definition a Killing form (dual to the Killing vector field $\xi$), which is equivalent to the first equation, and we see, under a re-scaling of the scalar curvature to $\s = n(n-1)$, that the Sasaki equation on an Einstein manifold insures that $\eta$ is moreover normal, that is, defines a $2$-form on tractors that is parallel with respect to $\nabla^{NC}$ (and this is where the conformal holonomy reduction comes from). Note, finally, applying the integrability conditions from $\curv^{\nabla^{NC}}$ for the normal conformal Killing form $\eta$, that we have $\W (X,Y)\eta = 0$, which is equivalent to $\xi \: \W = 0$. In other words, Einstein-Sasaki spaces also give a non-trivial example of Einstein spaces with degenerate Weyl tensor. We will discuss below the extent to which there are coverse statements to these facts, for which we will need some tools from contact geometry.

\section{Sasaki geometry and the $\k$-nullity distribution}

\v In general, an odd-dimensional manifold together with a 1-form, $(M^{2n+1},\eta)$ is called a {\it contact manifold} if the contact condition $\eta\wedge(d\eta)^n \neq 0$ holds. In this case, there exists a unique vector field $\xi$ satisfying
\begin{align*}
\xi \: d\eta = 0, \ \eta(\xi) = 1,
\end{align*}
\v called the {\it Reeb vector field} of the contact structure.\\

We are interested in contact structures arising in conformal Riemannian geometry, and so the metrics associated to contact structures are important. Following \cite{blair}, we define an {\it almost contact metric structure} as a set $(M,\eta,g,\phi)$ where, in addition to the contact structure, $g$ is a Riemannian metric and $\phi$ a (1,1)-tensor satisfying $\phi^2 = -Id + \eta\tens\xi$, such that
\begin{align*}
g(\phi X , \phi Y) = g(X,Y) - \eta(X) \eta(Y).
\end{align*}
\v If, in addition,
\begin{align}
g(X,\phi Y) = d\eta(X,Y),
\end{align}
\v then the structure is called {\it contact metric}. The Riemannian metric $g$ of a contact metric structure is called an {\it associated metric} to the contact structure $(M,\eta)$\\

Sasakian structures arise in contact geometry by imposing a further integrability condition on the contact metric structure, which is analogous to the step from almost K\"ahler manifolds to K\"ahler in symplectic geometry (cf. \cite{blair} p. 71). Thus, the first definition of a {\it Sasaki structure} is as a contact metric structure $(M,\eta,g,\phi)$ such that in addition the following integrability tensors vanish:
\begin{align*}
N^{(1)}(X,Y) &= [\phi,\phi](X,Y) + 2d\eta(X,Y)\xi,\\
N^{(2)}(X,Y) &= (\L_{\phi X}\eta)(Y) - (\L_{\phi Y}\eta)(X),\\
N^{(3)}(X) &= (\L_{\xi}\phi)(X),\\
N^{(4)}(X) &= (\L_{\xi}\eta)(X).
\end{align*}
\v And a standard result tells us that $N^{(1)}=0$ implies that the other tensors vanish as well.\\

There are other, equivalent, characterizations of Sasaki structures which occur in the literature, and it is useful for our purposes to review these. In general, for a contact metric structure, we have $N^{(2)}=N^{(4)}=0$. Further, the covariant derivative of $\xi$ satisfies
\begin{align*}
\nabla_X\xi = -\phi X - \frac{1}{2}\phi N^{(3)}X,
\end{align*}
\v and the tensor $N^{(3)}$ vanishes if and only if the Reeb vector field $\xi$ is a Killing field. A contact metric structure whose Reeb vector field is a Killing field is called {\it K-contact}, and for these, clearly $\phi = -\nabla\xi$.\\

It is a basic result of Sasakian geometry that the above definition of Sasaki structure is equivalent to the following: Let $(M,\eta,g,\phi)$ be an almost contact metric structure. $(M,\eta,g,\phi)$ is a Sasaki structure if, in addition:
\begin{align*}
(\nabla_X\phi)Y = g(X,Y)\xi - \eta(Y)X.
\end{align*}
 
From this, it is a straightforward computation to see that a K-contact metric structure is Sasakian if and only if:
\begin{align} \label{Sasaki curv}
\curv(X,Y)\xi = \eta(Y)X - \eta(X)Y,
\end{align}
\v This is often taken as the definition in the literature. The significance of this characterization for us is that it says, for an Einstein-Sasaki manifold, that the Reeb vector field $\xi \in Ker(\W)$.\\

Let us explain this point. The above equation says, in a special case ($\k=1$), that the Reeb vector field $\xi$ belongs to the so-called $\k$-nullity distribution (note that, for an almost contact metric structure, $\eta(Y) = g(\xi,Y)$). In general, on a pseudo-Riemannian manifold $(M,g)$ {\it the $\k$-nullity distribution} for a constant $\k$ is defined at a point $p \in M$ as:
\begin{align*}
N(\k)_p := \{Z \in T_pM : \curv(X,Y)Z = \k(g(Y,Z)X -g(X,Z)Y) \ \forall X,Y \in T_pM\}.
\end{align*}
\v $N(\k)_p$ is non-trivial for at most one constant $\k$, and enjoys some nice properties. Denoting by $n_{\k}(p)$ the rank of $N(\k)_p$ or {\it nullity index},  then although $n_{\k}$ is not constant, it is upper semi-continuous on $M$. Let $n_{\k,0}$ denote the minimum of the nullity index on $M$. Then $M_{n_{\k,0}} \subset M$ is open, where $M_{n_{\k,0}}$ denotes the subset where $N(\k)$ has rank $n_{\k,0}$. On $M_{n_{\k,0}}$, the $\k$-nullity distribution is integrable, i.e. a foliation, and totally geodesic. The leaves of the foliation are of constant curvature $\k$, and are complete if $(M,g)$ is complete. For proofs of these facts, see, e.g. \cite{gray70}, \cite{tanno78}.

\begin{Lemma} \label{Einstein nullity Lemma}
For an Einstein manifold $(M,g)$, the $\k$-nullity distribution $N(\k) = Ker(\W)$ and $\k = \frac {\s} {n(n-1)}$, where $\s$ is the scalar curvature.
\end{Lemma}

\v Proof. For the Schouten tensor $\h$ of an Einstein metric, the identity $\h = \frac {\s} {2n(n-1)}g$ holds. The Kulkarni-Nomizu product of the metric with itself acts on vectors as
\begin{align*}
g \star g (X,Y,U,V) = 2(g(X,U)g(Y,V)-g(X,V)g(Y,U)),
\end{align*}
\v Thus a vector $X$ lies in $Ker(\W = \curv - \h \star g)$ precisely when it lies in $N(\frac {\s} {n(n-1)})$.\\

Contact metric structures with a $\k$-nullity distribution (as well as more general distributions defined by curvature relations) have been studied quite a bit, especially in the case that the Reeb vector field $\xi$ lies in the $\k$-nullity distribution (cf. \cite{blair}, p. 105 for some references). To close this section, we cite one such result which will be applied in characterizing conformally Einstein-Sasaki spaces. This theorem of Tanno says that a contact metric structure on an Einstein manifold is Einstein-Sasaki precisely when the Reeb vector field lies in $Ker(\W)$:

\begin{Theorem} \label{Tanno Einstein Sasaki nullity Theorem}
(Tanno \cite{tanno88}): Let $(M,\eta,g,\phi)$ be a Riemannian Einstein manifold of dimension $2n+1\geq 5$ with contact metric structure. If $\xi$ belongs to the $\k$-nullity distribution, then $\k=1$ and $(M,\eta,g)$ is Sasakian.
\end{Theorem}

\section{Conformally Einstein-Sasaki spaces}

In the later sections, we will develop further conclusions from Lemma \ref{Einstein nullity Lemma} for the geometry of the distribution $Ker(\W)_{\vert M_k}$ for a conformally Einstein manifold with rank $k$ Weyl degeneracy, in order to prove the main result. But first, we show a simple application of this lemma and of Theorem \ref{Tanno Einstein Sasaki nullity Theorem} to the description of conformally Einstein-Sasaki spaces having rank 1 Weyl degeneracy. For a Riemannian manifold $(M,g)$ with rank 1 Weyl degeneracy, $Ker(\W)$ is a rank 1 distribution on an open and dense subset of $M$. Therefore, choosing a direction in $Ker(\W)$ and letting $\xi'$ at each point be the vector in this direction of unit-length, we get a uniqely determined unit-length vector field on this dense submanifold. Assume that $\xi'$ extends by taking limits to a smooth vector field on $M$. This clearly must be the case for a manifold $(M,g)$ conformally related to an Einstein-Sasaki space, and by continuity the extended vector field must also have constant unit-length. We denote this vector field by $\xi$, and its dual 1-form $\xi^*$ with respect to the metric $g$ we denote by $\eta$. The following proposition gives tensorial criteria for $(M,g)$ to be globally conformally Einstein-Sasaki.

\begin{Proposition}
Let $(M^{2n+1},g)$ be an oriented Riemannian manifold with rank 1 Weyl degeneracy . Then $(M,g)$ is conformally equivalent to an Einstein Sasaki space if and only if the following hold, with $\xi$ and $\eta = \xi^*$ as above:\\

\v 1. $\eta$ is contact ($\eta \wedge (d\eta)^{n} \neq 0$).\\

\v 2. For any positively-oriented orthonormal basis $\{\xi, E_1, \ldots, E_{2n}\}$, and $\Omega_{ij} = d\eta(E_i,E_j)$, the matrix $\frac {1} {\sqrt[n]{\vert det\Omega \vert}} \Omega$ is orthoganal.\\

\v 3. $\ric^o + (2-n)F_{(\nabla_g\sqrt[n]{\vert det\Omega \vert})} = 0$ (i.e., $\bar{g} = e^{(2\sqrt[n]{\vert det\Omega \vert})}g$ is an Einstein metric).
\end{Proposition}

\v Proof. First we note that the conditions of the proposition are conformally well-defined, i.e. they don't depend on which metric $g$ is chosen from the conformal class. Let $\bar{g} = e^{2\phi}g$ be a conformally related metric. First, $\eta$ depends conformally covariantly on $g$ (so the conformal structure gives us a conformal family of 1-forms):

\begin{align*}
\eta_{\bar{g}} = \bar{\eta} &= \bar{g}(\bar{\xi},.) = e^{2\phi}g(e^{-\phi}\xi,.) = e^{\phi}\eta = e^{\phi}\eta_g
\end{align*}

\v It follows then that 

\begin{align*}
\bar{\eta} \wedge (d\bar{\eta})^n = e^{\phi}\eta \wedge (de^{\phi}\wedge\eta + e^{\phi}d\eta)^n = e^{(n+1)\phi}(\eta \wedge (d\eta)^n),
\end{align*} 

\v and the contact condition (1) is conformally invariant.\\

Furthermore, $Ker(W)$ is in general a conformal invariant, and so in our case the orthogonal decomposition $T(M) = \R\xi \dsum \xi^{\perp}$ is conformally invariant. Let $\{E_1, \ldots, E_{2n}\}$ be a ONB, with respect to $g$, for $\xi^{\perp}$. Then $\{\bar{E_1}, \ldots, \bar{E_{2n}}\}$ likewise gives a ONB w.r.t $\bar{g}$ for $\bar{\xi}^{\perp}$, where $\bar{E_i} = e^{-\phi}E_i$. Using this conformal translation of bases for the invariant subspace $\xi^{\perp}$, we have:

\begin{align*}
d\bar{\eta}(\bar{E_i},\bar{E_j}) &= (de^{\phi}\wedge\eta) (e^{-\phi}E_i,e^{-\phi}E_j) + e^{\phi}d\eta (e^{-\phi}E_i,e^{-\phi}E_j)\\
&= 0 + e^{-\phi}d\eta(E_i,E_j)
\end{align*}

\v It follows that condition (2) is conformally invariant:

\begin{align*}
\frac {1} {\sqrt[n]{\vert det\bar{\Omega} \vert}} \bar{\Omega} &= \frac{1}{\sqrt[n]{\vert e^{-n\phi} det\Omega \vert}} e^{-\phi} \Omega \\
&= \frac{1}{\sqrt[n]{\vert det\Omega \vert}} \Omega.
\end{align*}

\v It is noted, as well, that condition (2) is not dependent on the choice of ONB (consistent with $\xi$), since any other such ONB is obtained by a rotation-matrix $P$, and for the resulting matrix $\Gamma$ corresponding to this new ONB, we have $\Gamma = P^{-1}\Omega P$. Thus $det\Gamma = det\Omega$ and $\frac {1} {\sqrt[n]{\vert det\Gamma \vert}} \Gamma = P^{-1} (\frac {1} {\sqrt[n]{\vert det\Omega \vert}} \Omega) P$ is also orthogonal.\\

And condition (3) is conformally invariant. Suppose $e^{(2\sqrt[n]{\vert det\Omega \vert})}g$ is an Einstein metric. Then

\begin{align*}
e^{(2\sqrt[n]{\vert det\bar{\Omega} \vert})}\bar{g} = e^{(2(-\phi)(\sqrt[n]{\vert det\Omega \vert}))}e^{2\phi}g = e^{(2\sqrt[n]{\vert det\Omega \vert})}g
\end{align*}

\v which is likewise Einstein, and thus (3) holds for $\bar{g}$.\\

Additionally, we remark that although in general the vector $\xi$ corresponding to a given metric $g$ in the conformal class is not the Reeb vector field determined by the contact form $\eta$ belonging to $g$, nevertheless this does hold for any Einstein metric in the conformal class. In general, $\xi$ will be the Reeb vector field of the contact form $\eta$ if and only if the integral curves of $\xi$ are geodesics, since:

\begin{align*}
2d\eta(\xi,X) &= \xi\eta(X) - X\eta(\xi) - \eta([\xi,X])\\
&= g(\nabla_{\xi}\xi,X) + g(\xi,\nabla_{\xi}X) - 0 - g(\xi,\nabla_{\xi}X - \nabla_X\xi)\\
&= g(\nabla_{\xi}\xi,X)
\end{align*}

\v Since $Ker(W)$ is totally geodesic, $\xi$ which is of constant length in $Ker(W)$ has geodesic integral curves.\\

\v ($\Leftarrow$): Assume that conditions (1)-(3) hold. Writing $\sigma = \sqrt[n]{\vert det\Omega \vert}$ we have the following facts for the metric $\bar{g} = e^{2\sigma}g$ in the conformal class:\\

\v (a) From condition (3), $\bar{g}$ is an Einstein metric.\\

\v (b) Thus, from the final remark above, $\bar{\xi} \in Ker(W)$ with $\bar{g}(\bar{\xi},\bar{\xi})=1$ is the Reeb vector field for the contact form $\bar{\eta}=\bar{g}(\bar{\xi},.)$.\\

\v (c) By condition (2), $(\bar{\Omega}_{ij}) = (d\bar{\eta}(\bar{E_i},\bar{E_j}))$ is orthogonal, for a $\bar{g}$-ONB $\{\bar{E_1},\ldots,\bar{E_{2n}}\}$ spanning the contact distribution $\bar{\xi}^{\perp}$.\\

We get an associated contact metric structure $(M,\bar{\eta},\bar{g},\phi)$ to $\bar{\eta}$ by defining $\phi(\bar{E_i}) = \bar{\Omega}_i^j\bar{E_j}$ on $\bar{\xi}^{\perp}$ and $\phi(\bar{\xi}) = 0$. Further, since the Reeb vector field $\bar{\xi}$ belongs to the $\k$-nullity distribution ($\k = \frac {\bar{\s}} {n(n-1)}$), Tanno's Theorem \ref{Tanno Einstein Sasaki nullity Theorem} tells us that $\bar{\s}=n(n-1)$ and $(M,\bar{\eta},\bar{g})$ is Einstein-Sasaki.\\

\v ($\Rightarrow$): Starting from an Einstein-Sasaki structure $(M,\xi,g)$, we have the contact form $\eta = g(\xi,.)$ and the endomorphism $\phi(X) = -\nabla_X\xi$, defining a contact metric structure. Then choosing a $\phi$-ONB $\{E_1,\ldots,E_{2n}\}$ for the contact distribution, the identity $g(X,\phi(Y)) = d\eta(X,Y)$ shows that $\Omega = d\eta(E_i,E_j)$ is the standard symplectic matrix, which is of course orthogonal. By characterisation \ref{Sasaki curv} of a Sasaki structure, the Killing vector field $\xi$ lies in $Ker(\W)$ (note that a Sasaki-Einstein space must have scalar curvature $\s = n(n-1)$, where $n$ is the dimension of the manifold). Thus conditions (1)-(3) are automatically satisfied, and the discussion above of their conformal invariance shows that they hold for all other metrics in the conformal class. $\qed$\\

\section{Conformally Einstein metrics with rank 1 Weyl degeneracy}

As the first step toward characterizing more general conformally Einstein metrics with degenerate Weyl tensor, we deal in this section with the most simple case of rank 1 degeneracy. We also prove in this section a number of results on the behavior of Riemannian foliations under conformal change of the metric, which hold for higher rank foliations and will be used in the next section. To begin, though, we make a simple observation on the 1-form $\eta$ arising in the rank 1 case, which allows us to divide the metrics into ''transversally integrable" and ''transversally non-integrable" classes. This is a natural next step from the contact condition $\eta \wedge (d\eta)^n \neq 0$, which can be seen as a maximal non-integrability condition on the distribution transverse to the Reeb vector field.\\

We saw in the previous section that the vector $\bar{\xi}$ is the Reeb vector field to the associated contact form $\bar{\eta}$ if and only if the integral curves defined by it are geodesics with respect to $\bar{g}$. In fact we can say more, namely that there is at most one metric in the conformal class for which this can occur. This is a result of the contact condition $\eta \wedge (d\eta)^n \neq 0$, and moreover the basic result holds even when this condition is weakened. Explicitly, we have:

\begin{Remark} \label{Rank 1 geodesic metric}
Let $(M,c)$ be a conformal Riemannian manifold with rank 1 Weyl degeneracy, and suppose $\eta \wedge d\eta$ is non-vanishing for the class of 1-forms $\eta$ defined in the previous section. Then there exists at most one metric $g \in c$ such that the integral curves of the associated vector field $\xi$ are geodesics (i.e., such that $Ker(\W)_{\vert M_1}$ is a totally geodesic foliation with respect to the metric $g$).
\end{Remark}

\v Proof. Suppose $g \in c$ is a metric in the conformal class for which the property holds. Then for a conformally related metric $\bar{g} = e^{2\phi}g$, we have $d\bar{\eta} = de^{-\phi} \wedge \eta + e^{-\phi}d\eta$. Thus, for any $\bar{E_i} \in \bar{\xi}^{\perp}$:

\begin{align*}
(\bar{\xi} \: d\bar{\eta})(\bar{E_i}) &= de^{-\phi}(\bar{\xi})\eta(\bar{E_i}) - e^{-\phi}\eta(\bar{\xi})(de^{-\phi})(\bar{E_i}) \\
&= 0 - e^{-2\phi}(de^{-\phi})(\bar{E_i})
\end{align*}

\v Thus, $\bar{\xi}$ has geodesic integral curves if and only if $(de^{-\phi})(\bar{E_i}) = 0$ for all $\bar{E_i} \in \bar{\xi}^{\perp}$, i.e. if and only if $\nabla \phi$ is proportional to $\bar{\xi}$. But precisely because of the condition $\eta \wedge d\eta \neq 0$, there are no non-trivial gradient vector fields satisfying this. For suppose $\nabla f_1 = f_2\xi$ for some smooth functions $f_1$ and $f_2$. Then $\eta = \frac{1}{f_2}df_1$ and $\eta \wedge d\eta = -\frac{1}{f_2}df_1 \wedge df_2 \wedge df_1 = 0$, a contradiction. \qed\\

On the other hand, the condition $\eta \wedge d\eta = 0$ is equivalent to the integrability of the complementary distribution $\xi^{\perp}$. In general, even if the 3-form $\eta \wedge d\eta$ is not identically zero on $M$, it can still have a non-trivial zero-set, and the zero-set could even contain an open subset in $M$. But we can consider the interior in $M$ of the set $\{\eta \wedge d\eta = 0\}$, and the open set $\{\eta \wedge d\eta \neq 0\}$. Then the disjoint union of these two open sets is dense in $M$. Thus, for the analysis of the local conformal geometry, it suffices to consider the following two distinct classes of manifolds:

\v 1. $(M,g)$ with $\eta \wedge d\eta = 0$, which we call {\it transversally integrable}. These have complementary foliations $\R\xi$ and $\xi^{\perp}$.\\

\v 2. $(M,g)$ with $\eta \wedge d\eta \neq 0$ on a dense subset of $M$, which we call {\it transversally non-integrable}. For these, there is a unique metric $\bar{g} \in [g]$ for which $\bar{\xi}$ has geodesic integral curves, and this metric coincides with the (unique!) Einstein metric in the conformal class.\\

This division into cases, we'll see, is important because the geometry in the first case is greatly simplified. In the second case, on the other hand, uniqueness of the ''geodesic metric" in the conformal class allows us to give explicit tensorial criteria for the metric to be locally conformally Einstein. First, though, we need some facts from foliation theory.\\

A {\it $k$-foliation} of an $n$-dimensional manifold $M$ is an integrable rank $k$ distribution. We'll denote the distribution by $\E$. When $M$ has in addition a Riemannian or conformal Riemannian metric $g$, we call $(M,\E,g)$ a {\it Riemannian foliation}. In this case, we have a canonical transversal distribution (in general, non-integrable) $\E^{\perp}$ given by the orthogonal complement and decomposition $TM = \E \dsum \E^{\perp}$.  We denote the projection maps by $\pi$ and $\pi^{\perp}$, respectively. The {\it first fundamental form} of a Riemannian foliation $\E$ is a map $h:\E \times \E \rightarrow \E^{\perp}$:

\begin{align*}
h(X,Y) := \pi^{\perp} \nabla_X Y
\end{align*}

\v Let $\xi_1, \ldots, \xi_k \in \E$ be local orthonormal vector fields spanning $\E$. The {\it mean curvature vector} of the foliation, which is independent of the choice of of spanning orthonormal vector fields, is defined to be the vector $H$:

\begin{align*}
H := \sum_{i=1}^k h(\xi_i,\xi_i)
\end{align*}

\v The geometrically important classes of Riemannian foliations are characterized using these objects. A foliation is called {\it totally geodesic} if $h = 0$. It is called {\it minimal} if $H=0$; and it is called {\it totally umbilical} if $h(X,Y) = g(X,Y)H$. In all these cases, the integrable sub-manifolds given by the leaves of the foliation have the corresponding geometric properties (totally geodesic, etc.). Obviously, a totally geodesic foliation is also minimal and umbilical. Further, it is well known that the converse is true. Only the property of being totally umbilical is preserved under conformal change of the Riemannian metric. We're interested in finding out how the other properties behave under conformal transformation.\\

\begin{Lemma}
Let $(M,\E,g)$ be a Riemannian $k$-foliation with mean curvature vector $H$. Given a conformal change of metric $\bar{g} = e^{2\phi}g$, the following identity holds for the mean curvature vector $\bar{H}$ of the associated Riemannian foliation $(M,\E,\bar{g})$:
\end{Lemma}
\begin{align} \label{Conformal mean curvature formula}
\bar{H} = e^{-2\phi}(H - k(\pi^{\perp} \nabla \phi)).
\end{align}

\v Proof. We have the following, well-known identity for the conformally changed Levi-Civita connection:

\begin{align*}
\bar{\nabla}_X Y = \nabla_X Y + d\phi(X)Y + d\phi(Y)X - g(X,Y) \nabla\phi.
\end{align*}

\v Furthermore, the distributions $\E$ and $\E^{\perp}$, as well as the associated projection maps, are unchanged by the conformal transformation. Thus:

\begin{align*}
\bar{H} =& \sum_{i=1}^k \pi^{\perp} \bar{\nabla}_{e^{-\phi}\xi_i} e^{-\phi}\xi_i\\
=& e^{-\phi} \sum_{i=1}^k \pi^{\perp} (\xi_i(e^{-\phi}) \xi_i + e^{-\phi} \nabla_{\xi_i} \xi_i + 2e^{-\phi} d\phi(\xi_i) \xi_i - e^{-\phi}g(\xi_i,\xi_i) \nabla \phi)\\
=& e^{-2\phi} \sum_{i=1}^k (h(\xi_i,\xi_i) - \pi^{\perp} \nabla \phi)\\
=& e^{-2\phi}(H - k(\pi^{\perp} \nabla \phi)) \qed\\
\end{align*}

In particular, if our foliation is minimal for some metric $g$ in the conformal class, we see that for all other metrics in the conformal class, the associated mean curvature vector must be the transversal projection of a gradient vector field. Explicitly, we have the direct corollary:\\

\begin{Corollary} \label{Mean curvature Corollary}
Let $(M,\E,g)$ be a Riemannian $k$-foliation. There exists a metric $\bar{g} \in [g]$ for which the conformally-related foliation $(M,\E,\bar{g})$ is minimal if and only if there exists a function $\phi$ such that $\pi^{\perp} \nabla \phi = \frac{1}{k}H$. If this holds, then a metric $\bar{g}$ for which the foliation is minimal is given by $\bar{g} = e^{2\phi}g$. We call $(M,\E,g)$ conformally minimal, and call $\bar{g}$ a minimal metric in the conformal class.
\end{Corollary}

The transformation formula (\ref{Conformal mean curvature formula}) also allows us to generalize the observation in Remark \ref{Rank 1 geodesic metric} to general Riemannian foliations. The following lemma gives a limit on the number of independent totally geodesic metrics (and thus, also, minimal metrics) which can exist in the conformal class.\\

\begin{Lemma} \label{Totally geodesic Lemma}
Let $(M,\E,g)$ be totally geodesic. Then for a second metric $\bar{g} = e^{2\phi}g$ in the conformal class, the conformally related foliation is totally geodesic if and only if $\nabla \phi \in \E$.
\end{Lemma}

\v Proof. Since $(M,\E,g)$ is totally geodesic and therefore totally umbilical, all the conformally related foliations must also be totally umbilical. Thus we have, for $X,Y \in \E$:

\begin{align*}
\bar{h} (X,Y) &= \bar{g}(X,Y) \bar{H} = e^{-2\phi}\bar{g}(X,Y)(H_g - k(\pi^{\perp} \nabla^g \phi))\\
&= -ke^{-2\phi}\bar{g}(X,Y)(\pi^{\perp} \nabla^g \phi)\\
\end{align*}

\v and the last line clearly vanishes if and only if $\nabla \phi \in \E$. \qed\\

To round out the material on foliations, we give here two Thoerems on totally geodesic foliations. These allow us, given limits on curvature tensors determined by a totally geodesic foliation, to make statements about the existence and properties of an integrable transversal distribution. The {\it mixed sectional curvature} of a Riemannian foliation $(M,\E,g)$ is defined to be $K(\pi X \wedge \pi^{\perp} Y)$ for vectors $X$ and $Y$, where $K(X \wedge Y)$ is the usual sectional curvature of $(M,g)$ for the plane determined by $\{X,Y\}$. Proofs of the Theorems can be found in  \cite{BejFar}, pp. 129-130.\\

\begin{Theorem} \label{Transversal integrability Theorem zero curvature}
Let $\E$ be a totally geodesic foliation on a Riemannian manifold $(M,g)$. If all mixed curvatures of $M$ at a point $x_0$ are positive, then the transversal distribution is not integrable.
\end{Theorem}

\begin{Theorem} \label{Transversal integrability Theorem positive curvature}
(Tanno 72). Let $\E$ be a totally geodesic foliation on a Riemannian manifold $(M,g)$. Suppose that all mixed sectional curvatures of $M$ vanish identically on $M$ and the transversal distribution $\E^{\perp}$ is integrable. Then the foliation defined by $\E^{\perp}$ is also totally geodesic.
\end{Theorem}

We are now ready to prove Theorem \ref{main} in the case of rank 1 Weyl degeneracy. As above, we let $\xi$ be a unit length vector field in $Ker(\W)$, and let $\eta$ be its dual 1-form, and $E_2,\ldots,E_n$ a local orthonormal frame spanning the transversal direction. The first proposition is for the transversally integrable case:\\

\begin{Proposition}
Let $(M,g)$ be a Riemannian manifold of rank 1 Weyl degeneracy which is transversally integrable. Suppose $(M,g)$ is conformally Einstein, with non-negative Einstein scaling. Then $(M,g)$ must be conformally Ricci-flat with decomposable conformal holonomy.
\end{Proposition}

\v Proof. Suppose $g$ is an Einstein metric. We've seen that for an Einstein manifold $(M,g)$, the distribution $Ker(\W)$ coincides with the $\k$-nullity distribution ($\k = \frac {\s} {n(n-1)}$). By definition, all mixed sectional curvatures of the $\k$-nullity distribution are constant and equal to $\k$. Since the $\k$-nullity distribution is totally geodesic and $(Ker(\W))^{\perp}$ is integrable in the transversally integrable case, we can directly apply the above theorems:  If $(M,g)$ has positive scalar curvature, then the integrability of $Ker(\W)^{\perp}$ contradicts Theorem \ref{Transversal integrability Theorem positive curvature}. If $(M,g)$ is Ricci-flat, then Theorem \ref{Transversal integrability Theorem zero curvature} implies that $(Ker(\W))^{\perp}$ is also totally geodesic. Thus, $(M,g)$ is locally isometric to a Riemannian product $h^1 + l^{(n-1)}$ (this product is even global if $(M,g)$ is complete), where $h^1$ is a flat metric on the one-dimensional leaves of $Ker(\W)$. It follows that $l^{(n-1)}$ is also Ricci-flat, and thus the conformal holonomy of $(M,[g])$ is decomposable. \qed \\

In the transversally non-integrable case, we use the uniqueness of the minimal metric in the conformal class to give local tensorial criteria. Namely, Lemma \ref{Mean curvature Corollary}, a local gradient field $V$ giving a transformation to a minimal metric must have the form $V = H + f\xi$ for some smooth function $f$. Applying local gradience, we have
$$ 0 = rot(V) = rot(H) + rot(f\xi) = rot(H) + df \wedge \eta + frot(\xi). $$

\v In particular, since $\eta \wedge d\eta$ is non-vanishing on a dense set in $M$, then on a dense set we always have local transversal vector fields $E_i, E_j \in \xi^{\perp}$ such that $d\eta(E_i,E_j) = rot(\xi)(E_i,E_j) \neq 0$. The local gradience of $V$ implies that our function $f$ must satisfy, for all such transversal vector fields: $f = - \frac {rot(H)(E_i,E_j)} {rot(\xi)(E_i,E_j)}$. Conversely, we could just define our function $f$ locally in this way. The uniqueness of a local gradient gradient field $V$ having the above form (which follows from the uniqueness of the minimal metric in the transversally non-integrable case) implies that $f$ does not depend on the choice of transversal vector fields. We thus get a smooth function on the dense subset $M_1$ and for $f$ thus defined we have:\\

\begin{Proposition}
Let $(M,g)$ be a Riemannian manifold with Weyl degeneracy of rank 1, $\xi, \eta$ as above, which is transversally non-integrable, i.e. $\eta \wedge d\eta \neq 0$ on a dense subset of $M$. Then $(M_1,g)$ is locally conformally Einstein if and only if $Ric^o + (2-n)F_{(H+f\xi)} = 0$, for $f$ defined through local equations as above.
\end{Proposition}

\v Proof. The $\k$-nullity distribution is totally geodesic (in particular, minimal) and $Ker(\W)$ is the $\k$-nullity distribution for an Einstein metric. Thus, $(M,g)$ is locally conformally Einstein only if $(M,Ker(\W),g)$ is locally conformally minimal. The previous considerations tell us that there is locally at most one minimal metric in the conformal class, and that the local gradient field which determines it has the form $H + f\xi$. Again, by the uniqueness of this minimal metric, it must also be a local conformal Einstein transformation, which is equivalent to the equation $Ric^o + (2-n)F_{H+f\xi} = 0$ being satisfied (see Section 2).\qed\\

Note that $(M_1,g)$ being globally conformally Einstein is equivalent to what was called a conformal almost Einstein structure on $(M,g)$ in Section 3 (i.e., it determines a parallel tractor). Alternatively, one could formulate results on all of $M$ by checking whether the limits of $f$ and $\xi$ on the complement of $M_1$ exist and satisfy the appropriate conditions. The only obstacle to a global result thus lies in the topology of the manifold and extending local gradient fields to global ones. We can conclude immediately:

\begin{Proposition}
Let $(M,g)$ be a simply connected Riemannian manifold with Weyl degeneracy of rank 1 which is transversally non-integrable. Then $(M,g)$ is conformally almost Einstein if and only if $Ric^o + (2-n)F_{(H+f\xi)} = 0$.
\end{Proposition}

\section{Rank $k$ Weyl degeneracy, $k \geq 2$}

The method used to solve the problem for rank 1 Weyl degeneracy must now be generalized and adapted for the higher rank degeneracy cases. We first divide the problem into "transversally non-integrable" cases (where we hope to have at most one possible solution) and "transversally integrable" cases (where we hope the geometry simplifies); Then we formulate tensorial criteria in the "transversally non-integrable" cases, giving necessary and sufficient conditions for a manifold to be conformally Einstein. We will see that new difficulties arise in both cases, but that it's possible to extend the method, thus giving a complete tensorial characterization of all conformally Einstein manifolds having rank $k$ Weyl degeneracy and non-negative Einstein scaling.\\

We will be dealing with the following general set-up for the remainder: $(M,\E,g)$ is a foliated Riemannian manifold, $\E$ is totally umbilical of rank $k$. We will write $\{ \xi_1, \ldots, \xi_k \}$ for local orthonormal vector fields spanning $\E$, and $\eta_1, \ldots, \eta_k$ will be their dual 1-forms, also given locally. Local orthonormal vector fields spanning the complement $\E^{\perp}$ will be denoted $E_j$, $j \in \{k+1, \ldots, n\}$. The second fundamental form $h$ and mean curvature vector field $H$ of the Riemannian foliation are as defined above. If we consider a conformal change of the metric $\bar{g} = e^{2\phi}g$, we get another Riemannian foliation $(M,\E,\bar{g})$, which we shall call conformally related. We'll denote the objects associated to this foliation by adding a bar: $\bar{\xi_i}, \bar{h}, \bar{H}$, etc.\\

From Lemma \ref{Totally geodesic Lemma}, we can introduce a class of transversally non-integrable metrics, in close analogy to the rank 1 case, which have at most one totally geodesic metric in the conformal class. Namely, we call $(M,g)$ {\it transversally non-integrable} if there exist no non-zero local gradients $\nabla \phi \in Ker(\W)$ and if $(Ker(\W))^{\perp}$ is non-integrable. Then for such metrics, there is (up to re-scaling by a constant) at most one metric in the conformal class which is Einstein. This follows from Lemma \ref{Totally geodesic Lemma}, since for an Einstein metric, $Ker(\W)$ is totally geodesic.\\

We will call $(M,g)$ {\it transversally integrable} if either there exists a non-zero local gradient vector field $\nabla \phi \in Ker(\W)$, or if $(Ker(\W))^{\perp}$ is integrable. Restricting to these two classes gives all possible local geomteries, as in the rank 1 case. As opposed to the rank 1 case, though, it of course does not follow immediately that the transversal distribution $Ker(\W)^{\perp}$ is integrable, and this is not necessarily the case. Further, it need not be the case that the rank 1 distribution defined by $\nabla \phi$ is totally geodesic with respect to any metric in the conformal class. Recall that the existence of a totally geodesic distribution with integrable transversal distribution (the assumption, along with conditions on the mixed sectional curvature, in Theorems \ref{Transversal integrability Theorem zero curvature} and \ref{Transversal integrability Theorem positive curvature}), is what allowed us to simplify the geometry for conformally Einstein metrics with non-negative Einstein scaling. For higher rank degeneracy, we can still do this, thanks to the following lemma on the $\k$-nullity distribution:\\

\begin{Lemma}
Let $(M,N(\k),g)$ be a Riemannian manifold with $\k$-nullity distribution of rank $n_{\k}$. Suppose there exists a 1-distribution contained in the $\k$-nullity distribution, $\E^1 = \{\xi_1\} \subset N(\k)$, such that $(\E^1)^{\perp}$ is involutive. Then locally there exists, for some $1 \leq k \leq n_{\k}$, a rank $k$ distribution $\E^k$, such that $\E^1 \subset \E^k \subset N(\k)$, $(\E^k)^{\perp}$ is involutive, and $\E^k$ is totally geodesic (and thus integrable as distribution).
\end{Lemma}

\v Proof. The proof is by induction on $k$:\\
\v ($k=1$): This would mean that $\E^1$ is totally geodesic, in which case we're done. If this isn't the case, then we can construct a rank 2 distribution $\E^2$, $\E^1 \subset \E^2 \subset N(\k)$, such that $(\E^2)^{\perp}$ is involutive: Since $\E^1$ is not totally geodesic, it must be that $\nabla_{\xi_1} \xi_1 \notin \E^1$. Then we define

\begin{align*}
\xi_2 := \frac {1} {\vert\vert \pi_1^{\perp}(\nabla_{\xi_1} \xi_1) \vert\vert} \pi_1^{\perp}(\nabla_{\xi_1} \xi_1).
\end{align*}

\v Here, $\pi_1$ denotes the projection associated to the distribution $\E^1$. We define $\E^2 := \{\xi_1,\xi_2\}$. Since $N(\k)$ is totally geodesic, $\xi_2 \in N(\k)$ and $\E^1 \subset \E^2 \subset N(\k)$. It remains to show that $(\E^2)^{\perp}$ is involutive. For this, since $(\E^1)^{\perp}$ is involutive, it clearly suffices to show that $g([E_i,E_j], \nabla_{\xi_1} \xi_1) = 0$ for all $E_i,E_j \in (\E^2)^{\perp}$.

\begin{align*}
- g([E_j,E_i], \nabla_{\xi_1} \xi_1) = g(\nabla_{\xi_1}([E_j,E_i]), \xi_1)\\
= \underbrace{g(\nabla_{\xi_1}\nabla_{E_j} E_i,\xi_1)} \ & \ \underbrace{-g(\nabla_{\xi_1}\nabla_{E_i} E_j,\xi_1)}\\
= g(\nabla_{E_j}\nabla_{\xi_1} E_i,\xi_1) - g(\nabla_{[E_j,\xi_1]} E_i,\xi_1) \ & \ - g(\nabla_{E_i}\nabla_{\xi_1} E_j,\xi_1) + g(\nabla_{[E_i,\xi_1]} E_j,\xi_1)\\
+g(\curv(\xi_1,E_i)E_j,\xi_1) \ & \ - g(\curv(\xi_1,E_j)E_i,\xi_1)\\
= \underbrace{g(\nabla_{E_j}\nabla_{\xi_1} E_i,\xi_1) + g(\nabla_{[E_i,\xi_1]} E_j,\xi_1)} \ & \ \underbrace{-g(\nabla_{E_i}\nabla_{\xi_1} E_j,\xi_1) - g(\nabla_{[E_j,\xi_1]} E_i,\xi_1)}\\
= g(\nabla_{\nabla_{\xi_1} E_i} E_j,\xi_1) + g(\nabla_{\nabla_{E_i} \xi_1} E_j,\xi_1) \ & \ - g(\nabla_{\nabla_{\xi_1} E_j} E_i,\xi_1) - g(\nabla_{\nabla_{E_j} \xi_1} E_i,\xi_1)\\ 
-g(\nabla_{\nabla_{\xi_1} E_i} E_j,\xi_1) \ & \ + g(\nabla_{\nabla_{\xi_1} E_j} E_i,\xi_1)\\
= g(\nabla_{\nabla_{E_i} \xi_1} E_j,\xi_1) - g(\nabla_{\nabla_{E_j} \xi_1} E_i,\xi_1)\\
= g(\nabla_{E_j}\nabla_{E_i} \xi_1,\xi_1) - g(\nabla_{E_i}\nabla_{E_j} \xi_1,\xi_1)\\
= g(\curv(E_j,E_i)\xi_1,\xi_1) + g(\nabla_{[E_j,E_i]} \xi_1,\xi_1)\\
= 0.
\end{align*}

\v Here we use the fact that $\xi_1 \in N(\k)$, to go from the 3rd line to the 4th, and then again in the final step of the calculation. Further, note that line 5 follows from line 4 because $[\nabla_{\xi_1} E_i,E_j] \in (\E^1)^{\perp}$ for all $E_i,E_j \in (\E^2)^{\perp}$. This holds, since $g(\nabla_{\xi_1} E_i,\xi_1) = -g(E_i,\nabla_{\xi_1} \xi_1) = 0$, based on the definition of $\E^2$, and since $(\E^1)^{\perp}$ is involutive. And line 7 follows from line 6 because $[\nabla_{E_i} \xi_1,E_j] \in (\E^1)^{\perp}$ for all $E_i,E_j \in (\E^2)^{\perp}$. This, in turn, follows from $g(\nabla_{E_i} \xi_1,\xi_1) = 0$ and the involutiveness of $(\E^1)^{\perp}$. Finally, note that from the involutiveness of $(\E^1)^{\perp}$, we also have $g([E_i,\xi_2],\xi_1) = 0$ for all $E_i \in (\E^2)^{\perp}$.\\

\v (Induction step): Assume that the above procedure has been carried out $(k-1)$ times, to get a series of distributions $\E^1 \subset \E^2 \subset \ldots \subset \E^k \subset N(\k)$, with $(\E^i)^{\perp}$ involutive, and $rank(\E^i) = i$ for all $i \in \{1,\ldots,k\}$. Let $\xi_i$ be a local unit-length vector field in $\E^i$ which is orthogonal to $\E^{(i-1)}$. Then from the involutiveness of $(\E^i)^{\perp}$ at each step, it follows that $g([E_i,\xi_r],\xi_s) = 0$ for all $1 \leq s \lneq r \leq k$ and for all $E_i \in (\E^s)^{\perp}$. If $\E^k$ is not totally geodesic, we will show how to construct $\E^{(k+1)} \subset N(\k)$ of rank $(k+1)$, containing $\E^k$, with $(\E^{(k+1)})^{\perp}$ involutive.\\

The key is the method for choosing a vector field to add to $\E^k$. Let $\nabla_{\xi_p} \xi_l \notin \E^k$ be the first such vector from the series: $\nabla_{\xi_1} \xi_1 \rightarrow \nabla_{\xi_2} \xi_1 \rightarrow \nabla_{\xi_1} \xi_2 \rightarrow \nabla_{\xi_2} \xi_2 \rightarrow \ldots \rightarrow \nabla_{\xi_{(k-1)}} \xi_{(k-1)} \rightarrow \nabla_{\xi_k} \xi_1 \rightarrow \ldots \rightarrow \nabla_{\xi_k} \xi_{(k-1)} \rightarrow \nabla_{\xi_1} \xi_k \rightarrow \ldots \rightarrow \nabla_{\xi_{(k-1)}} \xi_k \rightarrow \nabla_{\xi_k} \xi_k$. Then, analogous to the construction of $\E^2$, we define:

\begin{align*}
\xi_{(k+1)} := \frac {1} {\vert\vert \pi_k^{\perp}(\nabla_{\xi_p} \xi_l) \vert\vert} \pi_k^{\perp}(\nabla_{\xi_p} \xi_l),
\end{align*}

\v where $\pi_k$ is the projection associated to the distribution $\E^k$, and let $\E^{(k+1)} := \E^k \dsum \{ \xi_{(k+1)} \}$. To show the involutiveness of $(\E^{(k+1)})^{\perp}$, it suffices to show that $g([E_i,E_j],\nabla_{\xi_p} \xi_l) = 0$ for all $E_i,E_j \in (\E^{(k+1)})^{\perp}$. Note that, since $\xi_p, \xi_l \in N(\k)$, the proof from above for $k = 1$ goes through exactly the same, as long as we can show: (a) $[\nabla_{\xi_p} E_i,E_j] \in (\E^l)^{\perp}$; and (b) $[\nabla_{E_i} \xi_p,E_j] \in (\E^l)^{\perp}$ for all $E_i,E_j \in (\E^{(k+1)})^{\perp}$. This follows from our method of choosing $\nabla_{\xi_p} \xi_l$ from the above series. There are three cases to be checked: $l \lneq p$; $p \lneq l$; and $p=l$. We show (a) and (b) for the case $p \lneq l$, the other two cases being shown similarly:\\

Since $\nabla_{\xi_p} \xi_l$ is the first vector field from the series which is not contained in $\E^k$, in particular we must have:

\begin{align*}
\nabla_{\xi_l} \xi_i &\in \E^k, \ \forall i \lneq l,\\
\nabla_{\xi_i} \xi_l &\in \E^k, \ \forall i \lneq p,\\
\nabla_{\xi_i} \xi_j &\in \E^k, \ \forall i,j \lneq l.
\end{align*}

\v Let $E_j, E_k \in (\E^{(k+1)})^{\perp}$. Then we have $g(\nabla_{\xi_p} E_j,\xi_i) = - g(E_j,\nabla_{\xi_p} \xi_i) = 0$ for all $i \leq l$. So $\nabla_{\xi_p} E_j \in (\E^l)^{\perp}$, and (a) follows from involutiveness of $(\E^l)^{\perp}$.\\

From the involutiveness of $(\E^r)^{\perp}$ at each step, we also have: $g([E_j,\xi_s],\xi_r) = 0$ for all $1 \leq r \lneq s \leq k$. Thus, in particular we have:

\begin{align*}
g(\nabla_{E_j} \xi_s,\xi_p) &= g(\nabla_{\xi_s} E_j,\xi_p)\\
&= - g(E_j,\nabla_{\xi_s},\xi_p) = 0, \ \forall \ p \leq s \leq l.
\end{align*}

\v On the other hand, we also have:

\begin{align*}
g(\nabla_{E_j} \xi_p,\xi_r) &= g(\nabla_{\xi_p} E_j,\xi_r)\\
&= - g(E_j,\nabla_{\xi_p} \xi_r) = 0, \ \forall \ 1 \leq r \lneq p,
\end{align*}

\v and, since $\xi_p$ is unit-length, $g(\nabla_{E_j} \xi_p,\xi_p) = 0$. In sum, therefore, $g(\nabla_{E_j} \xi_p,\xi_i) = - g(\nabla_{E_j} \xi_i,\xi_p) = 0$ for all $i \leq l$, and $\nabla_{E_j} \xi_p \in (\E^l)^{\perp}$, from which (b) follows.\\

Finally, since $rank(N(\k)) = n_{\k}$ is finite, it is clear that either this process terminates at some step $k \lneq n_{\k}$, which means $\E^k$ is totally geodesic, or else $N(\k)^{\perp}$ is involutive. In any case, there must exist a totally geodesic (and thus involutive) distribution $\E^k$, $\E^1 \subset \E^k \subset N(\k)$, whose orthogonal complementary distribution is involutive. \qed \\

With this Lemma, we can now draw the following conclusion about the spaces we defined above as transversally integrable:

\begin{Proposition} \label{TransIntProp}
Let $(M,g)$ be a Riemannian manifold with rank $k$ Weyl degeneracy, $k \geq 1$, which is transversally integrable. Suppose $(M,g)$ is conformally Einstein, with non-negative Einstein scaling. Then $(M,g)$ is either locally conformally flat or conformally Ricci-flat with decomposable conformal holonomy.
\end{Proposition}

\v Proof. Let $\bar{g}$ be an Einstein metric in the conformal class. From the above Lemma, we have $\E^k \subset N(\frac {\bar{\s}} {n(n-1)})$ which is totally geodesic with $(\E^k)^{\perp}$ involutive. The mixed curvature of $\E^k$ is therefore constant and equal to $\frac {\bar{\s}} {n(n-1)}$. If $\bar{\s}$ is positive, then by Theorem 11, $(\E^k)^{\perp} = \{0\}$ and $(M,g)$ is conformally flat. If $\bar{\s} = 0$, then $(M,g)$ is conformally Ricci-flat. And by Theorem 12, $(\E^k)^{\perp}$ is also totally geodesic, and $\bar{g}$ is a local Riemannian product: $\bar{g} = \bar{g}_{\E^k} + \bar{g}_{(\E^k)^{\perp}}$. And since the leaves of $\E^k$ are flat, the leaves of the transversal foliation must be Ricci-flat, which shows that $(M,g)$ has decomposable conformal holonomy. \qed \\

Now that the transversally integrable case is taken care of, it remains to find tensorial criteria for a transversally non-integrable conformal manifold of arbitrary rank $k$ Weyl degeneracy to be (locally) conformally Einstein. The key fact which allows us to write down these conditions is the uniqueness of the metric in the conformal class with respect to which $Ker(\W)$ can be totally geodesic. As in the rank 1 case, uniqueness of the geodesic metric in the conformal class allows us first to define a unique local gradient vector field corresponding to this metric. The algorithm for defining this field through local equations is more complicated:

\begin{Remark} \label{Rank k algorithm} In this situation, there exists at most one vector field $V \in Ker(\W)$ such that that $\frac{1}{k}H + V$ is a local gradient. Moreover, a candidate $\overrightarrow{f}\overrightarrow{\xi} \in Ker(\W)$ can be determined by a procedure involving linear algebra and solving a system of first order differential equations.
\end{Remark}

\v Let $E_i, E_j \in Ker(\W)^{\perp}$. Then $\eta_p (E_i) = 0$ for the 1-forms dual to the $\xi_p \in Ker(\W)$. We thus have:

\begin{align*}
0 &= rot(\frac {1} {k} H + f_1 \xi_1 + \ldots + f_k \xi_k) (E_i,E_j)\\
&= \frac {1} {k} rot(H)(E_i,E_j) + f_1 rot(\xi_1)(E_i,E_j) + \ldots + f_k rot(\xi_k)(E_i,E_j)
\end{align*}

\v Taking pairs $(E_{i_1},E_{j_1}), \ldots, (E_{i_N},E_{j_N})$ for some $N$, with $i_k \lneq j_k$, we write $\Omega_p^r := rot(\xi_r)(E_{i_p},E_{j_p}) = d\eta_r(E_{i_p},E_{j_p})$. We thus have a $k \times N$ matrix

$$\Omega =
\left(\begin{array}{ccc} 
  \Omega_1^1 & \ldots & \Omega_1^k \\
  \Omega_2^1 & \ddots & \Omega_2^k \\
  \vdots & \vdots & \vdots \\
  \Omega_N^1 & \ldots & \Omega_N^k
\end{array}\right).$$

\v We need to determine uniquely a solution to the system $\Omega (f_{i=1}^N) = -(rot(H)_{i=1}^N)$. Note that, on a dense subset of $M$, the rank of this matrix (and similarly, those which appear later) must be constant. We assume, therefore, that we're at a point in this dense subset, for the following steps of the algorithm. As in the rank 1 case, for points where this assumption doesn't hold, it remains to see that the functions derived through this algorithm for the points in a neighborhood, when taking their limits, give smooth functions.\\

If $\Omega$ has maximal rank $k$, such a unique system is given by basic linear algebra. Suppose $\Omega$ has rank $r \lneq k$. Because $(M,g)$ is transversally non-integrable, $r$ must be non-zero. W.l.o.g., assume the matrix

$$\Omega(r) :=
\left(\begin{array}{ccc} 
  \Omega_1^1 & \ldots & \Omega_1^r \\
    \vdots & \vdots & \vdots \\
  \Omega_r^1 & \ldots & \Omega_r^r
\end{array}\right)$$

\v is non-singular. Then the fact that $\Omega$ has rank $r$ implies, for all $s \geq (r+1)$, that there exist functions $F^s_1, \ldots F^s_r$ such that $F^s_1 \Omega_j^1 + F^s_2 \Omega_j^2 + \ldots + F^s_r \Omega_j^r = \Omega_j^s$ for all $j \in \{1,\ldots,N\}$. We make the following change of coordinates $\{\xi_i\} \mapsto \{\xi_i'\}$ for $Ker(\W)$: 

\begin{align*}
\xi_l' &:= \xi_l + F^{(r+1)}_l \xi_{(r+1)} + F^{(r+2)}_l \xi_{(r+2)} + \ldots + F^k_l \xi_k, \ \forall 1 \leq l \leq r;\\
\xi_s' &:= \xi_s - F^s_1 \xi_1 - F^s_2 \xi_2 - \ldots - F^s_r \xi_r, \ \forall (r+1) \leq s \leq k.
\end{align*}

\v It can be checked that the $\xi_i'$ span $Ker(\W)$, and that $g(\xi_l',\xi_s') = 0$ for all $1 \leq l \leq r \lneq s \leq k$. After perhaps another change of basis, therefore, we can assume there are $\xi_i''$ orthonormal with $rot(\xi_s'')(E_i,E_j) = 0$ for all $(r+1) \leq s \leq k$. Then the linear relation above says that $rot(\xi_s)(E_i,E_j) = 0$ for all $E_i,E_j \in Ker(\W)^{\perp}$ and $(r+1) \leq s \leq k$. Thus we can rewrite the local gradient condition to get equations:

\begin{align*}
0 &= \frac {1} {k} rot(H)(E_i,E_j) + f_1'' rot(\xi_1'')(E_i,E_j) + \ldots + f_k'' rot(\xi_k'')(E_i,E_j)\\
&= \frac {1} {k} rot(H)(E_i,E_j) + f_1'' rot(\xi_1'')(E_i,E_j) + \ldots + f_r'' rot(\xi_r'')(E_i,E_j).
\end{align*} 

\v This gives a corresponding change $\Omega \mapsto \Omega''$, given at each point by multiplication by a non-singular $k \times k$ matrix. Thus $\Omega''$ has rank $r$, and we can find unique solutions for $f_1'',\ldots,f_r''$.

Next we can repeat an anlogous procedure, setting $H'' := \frac {1} {k} H + f_1''\xi_1'' + \ldots + f_r''\xi_r''$, and searching functions $f_{(r+1)}'', \ldots f_k''$ such that $rot(H'' + f_{(r+1)}\xi_{(r+1)}'' + \ldots + f_k''\xi_k'') = 0$. We can do this, determining some of the new functions, as long as $rot(\xi_s'')(X,Y) \neq 0$ for some vectors $X,Y \in \{\xi_{(r+1)}, \ldots, \xi_k''\}^{\perp}$. We continue in this manner until this is no longer the case. Then, either all functions have been determined and we're done, or, for some $p \geq (r+1)$, we have $(k-p+1)$ vectors in $Ker(\W)$ which we'll denote by $\xi_p^*, \ldots, \xi_k^*$ such that $rot(\xi_i^*)(X,Y) = 0$ for all $X,Y \in \{\xi_p^*, \ldots, \xi_k^*\}^{\perp}$. In this case, it can be checked that the only equations from the rotation free-equation which involve the still undetermined functions $f_p^*, \ldots, f_k^*$ are:

\begin{align*}
\xi_l^*(f_s^*) &- (f_p^* rot(\xi_p^*) + \ldots + f_k^* rot(\xi_k^*))(\xi_s^*,\xi_l^*)\\
&= \xi_s^*(f_l^*) + (f_1^* rot(\xi_1^*) + \ldots + f_{(p-1)}^* rot(\xi_{(p-1)}^*))(\xi_s^*,\xi_l^*);\\
\E_i(f_s^*) &- (f_p^* rot(\xi_p^*) + \ldots + f_k^* rot(\xi_k^*))(\xi_s^*,E_i)\\
&= rot(H)(\xi_s^*,E_i);\\
\xi_q^*(f_s^*) &- (f_p^* rot(\xi_p^*) + \ldots + f_k^* rot(\xi_k^*))(\xi_s^*,\xi_q^*)\\
&= 0;
\end{align*}

\v where $1 \leq l \leq (p-1) \lneq s,q \leq k$ and we write $f_1^*, \ldots, f_{(p-1)}^*$, respectively $\xi_1^*, \ldots \xi_{(p-1)}^*$, for the functions which have already been determined and the corrsponding orthonormal vectors. As usual, $E_i \in Ker(\W)^{\perp}$. Since the solution to this equation must be unique, the homogeneous part has only the trivial solution, which tells us that the matrix

$$\Omega^* :=
\left(\begin{array}{ccc} 
  rot(\xi_p^*)(\xi_p^*,\xi_1^*) & \ldots & rot(\xi_k^*)(\xi_p^*,\xi_1^*) \\
    \vdots & \vdots & \vdots \\
  rot(\xi_p^*)(\xi_k^*,\xi_{(p-1)}^*) & \ldots & rot(\xi_k^*)(\xi_k^*,\xi_{(p-1)}^*)\\
  rot(\xi_p^*)(\xi_p^*,E_{(k+1)}) & \ldots & rot(\xi_k^*)(\xi_p^*,E_{(k+1)}) \\
    \vdots & \vdots & \vdots \\
  rot(\xi_p^*)(\xi_k^*,E_n) & \ldots & rot(\xi_k^*)(\xi_k^*,E_n)\\
  rot(\xi_p^*)(\xi_p^*,\xi_p^*) & \ldots & rot(\xi_k^*)(\xi_p^*,\xi_p^*) \\
    \vdots & \vdots & \vdots \\
  rot(\xi_p^*)(\xi_k^*,\xi_k^*) & \ldots & rot(\xi_k^*)(\xi_k^*,\xi_k^*)\\
\end{array}\right)$$

\v has maximal rank $(k-p+1)$. Then taking a non-singular $(k-p+1) \times (k-p+1)$ sub-matrix of $\Omega^*$, we can solve the associated system of first order, linear differential equations associated to it, getting solutions $f_p^*, \ldots, f_k^*$.\\

\begin{Proposition} \label{Rank k transversally non-integrable}
Let $(M,g)$ be a Riemannian manifold of rank $k$ Weyl degeneracy which is transversally non-integrable. Then $(M_k,g)$ is locally conformally Einstein if and only if $Ric^o + (2-n)F_{\frac {1} {k} H + \overrightarrow{f}\overrightarrow{\xi}} = 0$, where $H$ is the mean curvature vector field for $Ker(\W)$ and $\overrightarrow{f}\overrightarrow{\xi}$ is the vector field defined in Remark \ref{Rank k algorithm}.
\end{Proposition}

\v Proof. ($\Rightarrow$): Let $(M_k,g)$ be locally conformally Einstein. Then there exists a local gradient vector field $\nabla \phi (= \nabla_g \phi)$ giving a conformal transformation to an Einstein metric $\bar{g} = e^{2\phi}g$. With respect to $\bar{g}$, $Ker(\W)$ is totally geodesic, in particular it is minimal. From Corollary 9, we must have: $\pi^{\perp} \nabla \phi = \frac {1} {k} H$, where $\pi = \pi_{Ker(\W)}$. Thus $\nabla \phi = \frac {1} {k} H + f_1 \xi_1 + \ldots + f_k \xi_k$ for some functions $f_1, \ldots, f_k$, and $rot(\frac {1} {k} H + f_1 \xi_1 + \ldots + f_k \xi_k) = 0$ since $\nabla \phi$ is a local gradient vector field. Further, since $(M,g)$ is transversally non-integrable, there can't be a second metric in the conformal class for which $Ker(\W)$ is totally geodesic, and so the $f_i$ are unique. \qed \\

As in the rank 1 case also, we have a global statement if $M$ is simply connected:

\begin{Proposition} \label{Rank k transversally non-integrable global}
Let $(M,g)$ be a Riemannian manifold of rank $k$ Weyl degeneracy which is transversally non-integrable and simply connected. Then $(M,g)$ is conformally almost Einstein if and only if $Ric^o + (2-n)F_{\frac {1} {k} H + \overrightarrow{f}\overrightarrow{\xi}} = 0$.
\end{Proposition}

\section{Applications, examples and further problems}

The algorithm given in Proposition \ref{Rank k transversally non-integrable}, at least in full generality, is rather cumbersome and not very elegant. It's therefore good to note that there are nice conditions which, while not sufficient, are necessary for a manifold with degenerate Weyl tensor to be conformally Einstein. The following restriction on the characteristic classes follows from the Theorem of Gray for the $\k$-nullity distribution (cf. \cite{gray70}).

\begin{Theorem}
Let $(M^n,g)$ be a Riemannian manifold of dimension $n$ with Weyl tensor of constant rank $k$ degeneracy, $k \geq 2$. If $(M,g)$ is conformally Ricci-flat, then $\chi(M) = 0$. If $(M,g)$ is conformally Einstein, then

\begin{align*}
P_i(M) = 0 \ \forall \ i \geq 1 + \frac {1} {4}(n-k)
\end{align*}

\v where $P_i(M) \in H^{4i}(M,\R)$ denotes the $i$th Pontryagin class.
\end{Theorem}

\v Proof. It clearly suffices to show this for $(M,g)$ a Ricci-flat (resp. any sort of Einstein) metric in the conformal class. Then $Ker(\W) = N(\k)$, the $\k$-nullity distribution for $\k = \frac {\s} {n(n-1)}$. Since this is of constant rank, we have a decomposition of the tangent bundle into vector bundles

\begin{align*}
T(M) = N(\k) \dsum N(\k)^{\perp}.
\end{align*}

\v Thus

\begin{align*}
\chi(M) = \chi(T(M)) = \chi(N(\k)) \chi(N(\k)^{\perp})
\end{align*}

\v and

\begin{align*}
P(M) = P(T(M)) = P(N(\k))P(N(\k)^{\perp}.
\end{align*}

\v Now in the Ricci-flat case, since $N(\k)$ is totally geodesic with flat leaves, $\chi(N(\k)) = 0$. In the general Einstein case, $N(\k)$ is totally geodesic with leaves of constant curvature. It follows that $P_i(N(\k)) = 0$ for all $i \geq 1$. Since $P_i(N(\k)) = 0$ for all $i \geq \frac{1}{4}(n-k)$, the result follows. \qed \\

The following Theorem restricts the number of independent Einstein metrics in the conformal class.

\begin{Theorem}
Let $(M,g)$ be an Einstein manifold with positive scalar curvature, and not of constant curvature. The only Einstein metrics $\bar{g} \in [g]$ in the conformal class are given by rescaling by a constant.
\end{Theorem}

\v Proof. Suppose $\bar{g} = \sigma g$ be another Einstein metric, $\sigma$ non-constant. Then $\sigma$ corresponds to a $\nabla^{NC}$-parallel tractor, i.e. $\sigma$ is a normal conformal Kiling 0-form. On the other hand, for an Einstein space, the exterior derivative of a normal conformal Killing p-form, if it's non-zero, is a normal conformal Killing (p+1)-form. Thus $d\sigma$ is a normal conformal Killing 1-form. In particular it follows from the integrability conditions that $\nabla \sigma \: \W = 0$. Thus, $Ker(\W)$ has rank $\geq 1$. But since both $g$ and $\sigma g$ are assumed to be Einstein, $(M,g)$ must be transversally integrable. This contradicts Proposition \ref{Rank k transversally non-integrable}, since $(M,g)$ has positive scalar curvature. \qed \\

\begin{Example}
A transversally integrable Einstein space with rank 1 Weyl degeneracy and negative scalar curvature.
\end{Example}

\v The above methods give no satisfactory solution for conformally Einstein spaces with negative Einstein scaling. We provide here an example of one such space, having degenerate Weyl tensor, which shows that this problem is not trivial. This also gives a non-trivial rank 1 example which is not an Einstein-Sasaki space. The mothod of construction is by warped products: Let $(M_*,g_*)$ be an $n$-dimensional Riemannian manifold and $f: I \rightarrow \R$ a nowhere vanishing smooth function on the interval. Then the Riemannian {\it warped product structure} defined by these data is

\begin{align*}
(I,dt^2) \times_f (M_*,g_*) = ((I \times M_*),(g = dt^2+f^2(t)g_*)).
\end{align*}

\v For a summary of the curvature and other geometric properties of the warped product, see \cite{KueRad}, pp. 6-7. Note, in particular, that as a special case we have the following

\begin{Lemma} \label{EinstWarp}
The warped product $(I,dt^2) \times_f (M_*,g_*)$ is an Einstein metric if and only if $g_*$ is an Einstein metric and $f'^2 + \rho f^2 = \rho_*$.
\end{Lemma}

\v Here, $\rho$ and $\rho_*$ denote the {\it normalized scalar curvatures} of the warped product metric resp. of $g_*$. We take $(M_*,g_*)$ to be a positively scaled Einstein space with non-degenerate Weyl tensor and $\rho_* = 1$. Then letting $f = sinh$, it can be directly checked that this gives $\rho = -1$ (formula (23) in \cite{KueRad}) and Lemma \ref{EinstWarp} implies that the warped product structure is Einstein with negative scalar curvature. Further, we have the following identities for the curvature tensor:

\begin{align}
\curv(X,Y)Z = \curv_*(X,Y)Z - \frac{f'^2}{f^2}\{g(Y,Z)X-g(X,Z)Y\} \label{one} \\
\curv(X,Y)\partial_t = 0 \label{two} \\
\curv(X,\partial_t)\partial_t = -\frac{f''}{f}X \label{three}
\end{align}

\v From identities (\ref{two}) and (\ref{three}) it follows immediately that the vector field $\partial_t$ lies in the $\k$-nullity distribution ($\k = \rho = -1$). Moreover, this warped product structure does not have higher rank Weyl degeneracy, since that would imply, by identitiy (\ref{one}), that

\begin{align}
\curv_*(X,Y)Z = (\frac{f'^2}{f^2} - 1)\{g(Y,Z)X - g(X,Z)Y\} \label{four}
\end{align}

\v for some $Z \in \partial_t^{\perp}$ and for all vector fields $X,Y$. But considering a level set submanifold $\{t = c_0\}$ and vector fields $X,Y$ tangent to the level set, identity (\ref{four}) implies

\begin{align*}
\curv_*(X,Y)Z = (\frac{f'^2(c_0)}{f^2(c_0)} - 1)f^2(c_0)\{g_*(Y,Z)X - g_*(X,Z)Y\},
\end{align*}

\v which is excluded, since $(M_*,g_*)$ was taken to be Einstein with non-degenerate Weyl tensor, and thus has no non-trivial $\k$-nullity distribution. \qed \\

This example shows that conformally Einstein spaces exist with negative Einstein scaling and degenerate Weyl tensor. However, this example does not show that a result as in Proposition \ref{TransIntProp} can not be extended to transversally integrable conformally Einstein spaces with negative Einstein scaling. Indeed, in our example, conformally rescaling the warped product metric by the warping factor gives a second Einstein metric in the conformal class. Thus, we see that this metric has decomposable conformal holonomy. An open question (as far as we know) is whether negative Einstein metrics can be constructed which are degenerate and have non-decomposable conformal holonomy. This would also be an interesting question to answer, because such examples would provide a case where no (normal) conformal Killing forms are possible - as opposed to the examples we've seen, such as Sasaki-Einstein and conformally decomposable spaces. Alternatively, we plan to look into a method for characterizing transversally integrable conformal spaces with negative Einstein scaling. Further steps in our investigation include the extension of these results to the Lorentzian setting, which would include the important case of conformally Einstein Feffermann spaces.

\end{document}